\documentclass{amsart}
\usepackage{amsfonts}
\usepackage{graphicx}

\newtheorem{theorem}{Theorem}[section]
\newtheorem*{theorem A}{Theorem A}
\newtheorem*{theorem B}{N\"olker's Theorem}

\newtheorem{proposition}{Proposition}[section]
\newtheorem{corollary}{Corollary}[section]

\theoremstyle{remark}

\theoremstyle{remark}

\theoremstyle{definition}

\newtheorem{example}{Example}[section]
\numberwithin{equation}{section}
\def\({\left ( }
\def\){\right )}
\def\<{\left < }
\def\>{\right >}

\input{tcilatex}

\begin{document}
\title{Constant curvature surfaces in a pseudo-isotropic space}
\author{Muhittin Evren Aydin}
\address{Department of Mathematics, Faculty of Science, Firat University,
Elazig, 23119, Turkey}
\email{meaydin@firat.edu.tr}
\thanks{}

\begin{abstract}
In this study, we deal with the local structure of curves and surfaces
immersed in a pseudo-isotropic space $\mathbb{I}_{p}^{3}$ that is a
particular Cayley-Klein space. We provide the formulas of curvature, torsion
and Frenet trihedron in order for spacelike and timelike curves. The causal
character of all admissible surfaces in $\mathbb{I}_{p}^{3}$ has to be
timelike or lightlike up to its absolute. We introduce the formulas of
Gaussian and mean curvature for timelike surfaces in $\mathbb{I}_{p}^{3}$.
As applications, we describe the surfaces of revolution which are the orbits of
a plane curve under a hyperbolic rotation with constant Gaussian
and mean curvature.
\end{abstract}

\maketitle

\section{Introduction and preliminaries}

Let $P\left( \mathbb{R}^{3}\right) $ be the projective 3-space and $\left(
x_{0}:x_{1}:x_{2}:x_{3}\right) $ the homogenous coordinates. By a \textit{%
quadric}, we mean a subset of points of $P\left( \mathbb{R}^{3}\right) $
described as zeros of a quadratic form associated with a non-zero symmetric
bilinear form of $P\left( \mathbb{R}^{3}\right) .$

The\textit{\ Cayley-Klein 3-spaces }can be defined in $P\left( \mathbb{R}%
^{3}\right) $ with an \textit{absolute figure}, namely a sequence of
quadrics and subspaces of $P\left( \mathbb{R}^{3}\right) $, see \cite%
{12,15,25,28}. We are interested in a particular Cayley-Klein space, the 
\textit{pseudo-isotropic 3-space} $\mathbb{I}_{p}^{3}$. Its \textit{absolute 
}is composed of the quadruple $\left\{ \omega ,f_{1},f_{2},F\right\} ,$
where $\omega $ is the plane at infinity, $f_{1},f_{2}$ two real lines in $%
\omega $, $F$ the intersection of $f_{1}$ and $f_{2}.$ In coordinate form,
these arguments are given by%
\begin{equation*}
\omega :x_{0}=0,\text{ }f_{1}:x_{0}=x_{1}=0,\text{ }f_{2}:x_{0}=x_{2}=0,%
\text{ }F\left( 0:0:0:1\right) .
\end{equation*}%
For further details, see \cite{8,13,17,18}.

We deal with an affine model of $\mathbb{I}_{p}^{3}$ via the coordinates $%
\left( x=\frac{x_{1}}{x_{0}},y=\frac{x_{2}}{x_{0}},z=\frac{x_{3}}{x_{0}}%
\right) ,$ $x_{0}\neq 0.$ The \textit{group of pseudo-isotropic motions} is
a six-parameter group given by%
\begin{equation}
\left( x,y,z\right) \longmapsto \left( x^{\prime },y^{\prime },z^{\prime
}\right) :\left\{ 
\begin{array}{l}
x^{\prime }=a+qx, \\ 
y^{\prime }=b+\frac{1}{q}y\text{ }\left( q\neq 0\right) , \\ 
z^{\prime }=c+dx+ey+z,%
\end{array}%
\right.  \tag{2.1}
\end{equation}%
where $a,b,c,d,e,q\in \mathbb{R}.$ The \textit{pseudo-isotropic metric }is
introduced by the absolute, i.e. $ds^{2}=dx^{2}-dy^{2}$. Note that this
metric can be also considered as $ds^{2}=dxdy$ by standing $x=(x+y)/2$, $%
x=(x-y)/2$.

The investigation of curves and surfaces in 3-spaces is a classical field of
study in differential geometry. In spite of the fact that the cyclides in $%
\mathbb{I}_{p}^{3}$, i.e. algebraic surfaces of order 4, have been studied
for many years; as far as we know, the local structure of curves and
surfaces in $\mathbb{I}_{p}^{3}$ have not been established.

Indeed, we found motivation for this paper in B. Divjak's works (\cite%
{8,9,20}), in which the author introduced the differential geometry of
curves and surfaces in the pseudo-Galilean space as generalizing that of the
Galilean space. Intending a similar approach for the isotropic geometry (for
details, see \cite{1,2,3,5,11,14,22,23,26,27}), we are interested in the
local theory of curves and surfaces in $\mathbb{I}_{p}^{3}.$

The fact that the pseudo-isotropic metric is indefinite requires to
introduce some basic notions (e.g. the causal character, the pseudo-angle,
etc.) in $\mathbb{I}_{p}^{3}$ from the semi-Riemannian geometry (see Section
2). For detailed properties of such a geometry see \cite{6,16,24}.

In Section 3, it is suprisingly observed that each lightlike curve in $%
\mathbb{I}_{p}^{3}$ lies in the isotropic plane of the form $x\pm y=c,$ $%
c\in \mathbb{R}.$ As the local structures of the non-lightlike curves, the
formulas in $\mathbb{I}_{p}^{3}$ analogous to the famous Frenet's formulas
were given.

We get in Section 4 that each immersed admissible surface in $\mathbb{I}%
_{p}^{3}$ is timelike or lightlike. The formulas of the Gaussian and the
mean curvatures for timelike surfaces are also introduced.

As several applications, in Section 5, we study and classify the surfaces of revolution, imposing some natural curvature conditions.

\section{Basics in the sense of pseudo-isotropic geometry}

The\textit{\ pseudo-isotropic scalar product }between two vectors $u=\left(
u_{1},u_{2},u_{3}\right) $ and $v=\left( v_{1},v_{2},v_{3}\right) \in 
\mathbb{I}_{p}^{3}$\textit{\ }can be defined as%
\begin{equation}
\left\langle u,v\right\rangle =\left\{ 
\begin{array}{ll}
u_{3}v_{3}, & \text{if }u_{1}=u_{2}=v_{1}=v_{2}=0\text{,} \\ 
u_{1}v_{1}-u_{2}v_{2}, & \text{otherwise.}%
\end{array}%
\right.  \notag
\end{equation}

A line is said to be \textit{isotropic} (resp. \textit{non-isotropic}) if
its point at infinity is (resp. no) the absolute point $F$. Moreover, a
plane is said to be \textit{isotropic} (resp. \textit{non-isotropic}) if its
line at infinity containes (resp. does not) the absolute point $F$. In the
affine model of $\mathbb{I}_{p}^{3}$, the isotropic lines and planes are
parallel to the $z-$axis. In the non-isotropic planes, the Lorentzian metric
is basically used.

Let us consider the projection onto $xy-$plane given by 
\begin{equation*}
u=\left( u_{1},u_{2},u_{3}\right) \longmapsto \tilde{u}=\left(u_{1},u_{2},0%
\right) ,
\end{equation*}
usually called \textit{top view}. A nonzero vector $u$ is said to be \textit{%
isotropic} (resp. \textit{non-isotropic}) if $\tilde{u}=0$ (resp. $\tilde{u}
\neq 0$). The zero vector is assumed to be non-isotropic. A non-isotropic
vector $u\in \mathbb{I}_{p}^{3}$ is respectively called \textit{spacelike}, 
\textit{timelike} and \textit{lightlike }(or\textit{\ null}) if $%
\left\langle u,u\right\rangle >0$ or $u=0,$ $\left\langle u,u\right\rangle
<0 $ and $\left\langle u,u\right\rangle =0$ $\left( u\neq 0\right) .$

The set of all lightlike vectors of $\mathbb{I}_{p}^{3}$ is called \textit{%
lightlike cone}, i.e.,%
\begin{equation*}
\Lambda =\left\{ \left. \left( u_{1},u_{2},u_{3}\right) \in \mathbb{I}%
_{p}^{3}\right\vert u_{1}^{2}-u_{2}^{2}=0\right\} -\left\{ 0\in \mathbb{I}%
_{p}^{3}\right\} .
\end{equation*}%
Denote $\mathcal{T}$ the set of all timelike vectors in $\mathbb{I}_{p}^{3}.$
For some $u\in \mathcal{T}$, the set given by%
\begin{equation*}
\mathcal{C}\left( u\right) =\left\{ v\in \mathcal{T}:\left\langle
u,v\right\rangle <0\right\}
\end{equation*}%
is called the \textit{timelike cone} of $\mathbb{I}_{p}^{3}$ containing $u.$

The \textit{pseudo-isotropic} \textit{angle} of two timelike non-isotropic
vectors $u,v\in \mathbb{I}_{p}^{3}$ lying in the same timelike-cone is
defined as the Lorentzian angle between $\tilde{u}$ and $\tilde{v}$, i.e. 
\begin{equation*}
\left\langle u,v\right\rangle =-\sqrt{-\left\langle u,u\right\rangle }\sqrt{%
-\left\langle v,v\right\rangle } \cosh \phi .
\end{equation*}%
Note that all isotropic vectors are isotropically orthonogal to
non-isotropic ones. Further, two non-isotropic vectors $u,v$ in $\mathbb{I}%
_{p}^{3}$ are orthonogal if $\left\langle u,v\right\rangle =0.$

\section{Spacelike and timelike curves in $\mathbb{I}_{p}^{3}$}

Let $\alpha \left( s\right) =\left( x\left( s\right) ,y\left( s\right)
,z\left( s\right) \right) $ be a regular curve in $\mathbb{I}_{p}^{3},$ i.e. 
$\alpha ^{\prime }\left( s\right) =\dfrac{d\alpha }{ds}\neq 0$ for all $s.$
Then it is said to be \textit{admissible} if $\alpha \left( s\right) $ has
no isotropic osculating plane. An admissible curve $\alpha \left( s\right) $ in $\mathbb{I}_{p}^{3}$ is
said to be \textit{spacelike} (resp. \textit{timelike}, \textit{lightlike})%
\textit{\ }if $\alpha ^{\prime }\left( s\right) $ is spacelike (resp.
timelike, lightlike) for all $s.$ An easy compute shows that all lightlike
curves lie in the isotropic plane of the form $x\pm y=c,$ $c\in \mathbb{R}.$

Henceforth, we consider only spacelike and timelike admissible curves.

Now let $\alpha =\alpha \left( s\right) $ be a spacelike curve in $\mathbb{I}%
_{p}^{3}$ parameterized by arc-length. Then we have 
\begin{equation}
\left\langle \alpha ^{\prime },\alpha ^{\prime }\right\rangle =\left(
x^{\prime }\right) ^{2}-\left( y^{\prime }\right) ^{2}=1  \tag{3.1}
\end{equation}%
and taking derivative of $\left( 3.1\right) $ gives%
\begin{equation}
x^{\prime }x^{\prime \prime }-y^{\prime }y^{\prime \prime }=0.  \tag{3.2}
\end{equation}%
Denote $T=\alpha ^{\prime }$ and call it \textit{tangent vector field.}
Since $T^{\prime }=\alpha ^{\prime \prime }$ is timelike in $\mathbb{I}%
_{p}^{3}$ we can define the following%
\begin{equation*}
\kappa =\sqrt{\left( y^{\prime \prime }\right) ^{2}-\left( x^{\prime \prime
}\right) ^{2}},
\end{equation*}%
called \textit{curvature} of $\alpha .$ Using $\left( 3.2\right) ,$ we get%
\begin{equation}
\kappa =\frac{y^{\prime \prime }}{x^{\prime }}\text{ or }\kappa =\frac{%
x^{\prime \prime }}{y^{\prime }}.  \tag{3.3}
\end{equation}%
Considering $\left( 3.1\right) $ and $\left( 3.2\right) $ into $\left(
3.3\right) $ we find%
\begin{equation}
\kappa =\det \left( \tilde{\alpha}^{\prime },\tilde{\alpha}^{\prime \prime
}\right) .  \tag{3.4}
\end{equation}

Define the \textit{normal vector field} and\textit{\ torsion} of $\alpha $
respectively as%
\begin{equation}
N =\frac{1}{\kappa }T^{\prime }\text{ and }\tau =\frac{\det \left( \alpha
^{\prime },\alpha ^{\prime \prime },\alpha ^{\prime \prime \prime }\right) }{%
\kappa ^{2}}, \text{ } \kappa \neq 0.  \tag{3.5}
\end{equation}

Since $B=\left( 0,0,1\right) $ is isotropically orthogonal to $T$ and $N,$
we can take it as the \textit{binormal vector field of }$\alpha .$

From $\left(3.5\right)$ we have%
\begin{equation}
N^{\prime } =\left( \frac{1}{\kappa }\right) ^{\prime }\left( x^{\prime
\prime } ,y^{\prime \prime } ,z^{\prime \prime } \right) +\frac{1}{\kappa }%
\left( x^{\prime \prime \prime } ,y^{\prime \prime \prime } ,z^{\prime
\prime \prime } \right) .  \tag{3.6}
\end{equation}%
Put $N^{\prime}=\left( n_{1},n_{2} ,n_{3} \right) $. Hence we write%
\begin{equation}
n_{1}=\left( \frac{1}{\kappa }\right) ^{\prime }x^{\prime \prime }+\frac{1}{%
\kappa }x^{\prime \prime \prime }.  \tag{3.7}
\end{equation}%
Using $\left( 3.4\right) $ into $\left( 3.7\right) $ yields%
\begin{equation}
n_{1}=-\frac{x^{\prime }}{\kappa ^{2}}\left( x^{\prime \prime }y^{\prime
\prime \prime }-x^{\prime \prime \prime }y^{\prime \prime }\right) . 
\tag{3.8}
\end{equation}%
By taking derivative of $\left( 3.2\right) $ and considering into $\left(
3.8\right) $ we obtain 
\begin{equation}
n_{1}=-\kappa x^{\prime }.  \tag{3.9}
\end{equation}%
Similar computations gives%
\begin{equation}
n_{2}=-\kappa y^{\prime }.  \tag{3.10}
\end{equation}%
For the third component of $N^{\prime}$, we have%
\begin{equation*}
n_{3}=\left( \frac{1}{\kappa }\right) ^{\prime }z^{\prime \prime }+\frac{1}{%
\kappa }z^{\prime \prime \prime }.
\end{equation*}%
It follows from $\left( 3.4\right) $ that%
\begin{equation}
n_{3}=\frac{1}{\kappa ^{2}}\left\{ -\left( x^{\prime }y^{\prime \prime
\prime }-x^{\prime \prime \prime }y^{\prime }\right) z^{\prime \prime
}+\left( x^{\prime }y^{\prime \prime }-x^{\prime \prime }y^{\prime }\right)
z^{\prime \prime \prime }\right\} .  \tag{3.11}
\end{equation}%
By adding and substracting $\left( x^{\prime \prime }y^{\prime \prime \prime
}-y^{\prime \prime }x^{\prime \prime \prime }\right) z^{\prime }$ in $\left(
3.11\right) $ we conclude%
\begin{equation}
n_{3}=\frac{1}{\kappa ^{2}}\left\{ \det \left( \alpha ^{\prime },\alpha
^{\prime \prime },\alpha ^{\prime \prime \prime }\right) -\left( x^{\prime
\prime }y^{\prime \prime \prime }-x^{\prime \prime \prime }y^{\prime \prime
}\right) z^{\prime }\right\} .  \tag{3.12}
\end{equation}%
Taking derivative of $\left( 3.2\right) $ and considering into $\left(
3.12\right) $ implies%
\begin{equation}
n_{3}=\tau -\kappa z^{\prime }.  \tag{3.13}
\end{equation}%
$\left( 3.9\right) ,$ $\left( 3.10\right) $ and $\left( 3.13\right) $ yield
that $N^{\prime }=-\kappa T+\tau B$. Thus we obtain the formulas analogous
to these of Frenet as follows%
\begin{equation}
\left\{ 
\begin{array}{l}
T^{\prime }=\kappa N \\ 
N^{\prime }=-\kappa T+\tau B \\ 
B^{\prime }=0.%
\end{array}%
\right.  \notag
\end{equation}

By similar arguments, we can find the derivative formulas of the vector
fields $T,N,B$ for a timelike curve in $\mathbb{I}_{p}^{3}$ as

\begin{equation}
\left\{ 
\begin{array}{l}
T^{\prime }=\kappa N \\ 
N^{\prime }=\kappa T-\tau B \\ 
B^{\prime }=0,%
\end{array}%
\right.  \notag
\end{equation}
where $\kappa =\sqrt{\left( x^{\prime \prime }\right) ^{2}-\left( y^{\prime
\prime }\right) ^{2}}$ and $\tau =\dfrac{\det \left( \alpha ^{\prime
},\alpha ^{\prime \prime },\alpha ^{\prime \prime \prime }\right) }{\kappa
^{2}}.$

\begin{example}
Consider a hyperbolic cylindrical curve in $\mathbb{I}_{p}^{3}$ given by
(see \cite{7})%
\begin{equation}
\alpha \left( s\right) =\left( \cosh s,\sinh s,z\left( s\right) \right) . 
\tag{3.14}
\end{equation}%
This is a timelike curve of arc-length in $\mathbb{I}_{p}^{3}$ with $\kappa
\left( s\right) =1$ and 
\begin{equation}
\tau \left( s\right) =z^{\prime }\left( s\right) -z^{\prime \prime \prime
}\left( s\right) .  \tag{3.15}
\end{equation}%
If $\alpha $ has constant torsion $\tau _{0},$ then by solving $\left(
3.15\right) $ we find%
\begin{equation*}
z\left( s\right) =\tau _{0}s+c_{1}e^{s}-c_{2}e^{-s}+c_{3},\text{ }%
c_{1},c_{2},c_{3}\in \mathbb{R},
\end{equation*}
which gives the elemantary result:
\end{example}

\begin{proposition}
Let $\alpha $ be a hyperbolic cylindrical curve in $\mathbb{I}_{p}^{3}$ with
constant torsion $\tau _{0}.$ Then it is of the form%
\begin{equation*}
\alpha \left( s\right) =\left( \cosh s,\sinh s,\tau
_{0}s+c_{1}e^{s}-c_{2}e^{-s}+c_{3}\right) ,
\end{equation*}%
where $c_{1},c_{2},c_{3}\in \mathbb{R}.$
\end{proposition}

\begin{figure}[ht]
\begin{center}
\includegraphics[scale=0.5]{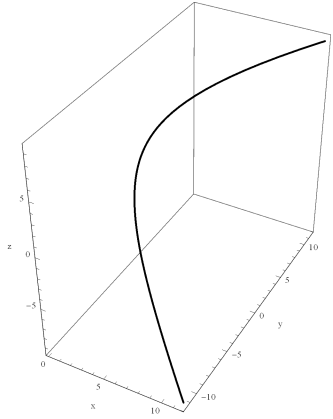}
\end{center}
\caption{A hyperbolic cylindrical curve with $\protect\tau _{0}=3$, $s \in (-%
\protect\pi , \protect\pi )$.}
\end{figure}

\section{Timelike surfaces in $\mathbb{I}_{p}^{3}$}

Let $M$ be a surface immersed in $\mathbb{I}_{p}^{3}$ without isotropic
tangent planes. Then we call such a surface \textit{admissible}. Let $T_{x}M$
be a non-isotropic tangent plane at a point $x \in M$. An admissible surface 
$M$ is said to be \textit{timelike} (resp. \textit{lightlike}) if the
induced metric $g$ in $T_{x}M$ for each $x \in M$ from $\mathbb{I}_{p}^{3}$
is non-degenerate of index 1 (resp. degenerate).

Henceforth, we will not consider lightlike surfaces.

Assume that $M$ has a local parameterization in $\mathbb{I}_{p}^{3}$ as
follows:%
\begin{equation*}
r:D\subseteq \mathbb{R}^{2}\longrightarrow \mathbb{I}_{p}^{3}:\text{ }\left(
u_{1},u_{2}\right) \longmapsto \left( x\left( u_{1},u_{2}\right) ,y\left(
u_{1},u_{2}\right) ,z\left( u_{1},u_{2}\right) \right)
\end{equation*}%
for smooth real-valued functions $x,y,z$ on a domain $D\subseteq \mathbb{R}%
^{2}.$ Denote $\left( g_{ij}\right) $ the matrical expression of $g$ with
respect to the basis $\left\{ r_{u_{1}},r_{u_{2}}\right\} .$ Then we have 
\begin{equation}
g_{ij}=\left\langle r_{u_{i}},r_{u_{j}}\right\rangle ,\text{ }r_{u_{i}}=%
\frac{\partial r}{\partial u_{i}},\text{ }i,j=1,2.  \notag
\end{equation}%
It is easy to see that 
\begin{equation*}
\det
\left(g_{ij}\right)=-\left(x_{u_{1}}y_{u_{2}}-x_{u_{2}}y_{u_{1}}\right)^{2}.
\end{equation*}
The unit normal vector field of $M$ is the isotropic vector $\xi=\left(
0,0,1\right) $ since it is isotropically orthogonal to the tangent plane of $%
M$.

For the second fundamental form of $M$, we follow the similar way with Sachs
(see \cite{20}, p. 155). Let $r\left( s\right) $ be an arc-length curve on $%
M $ and $T$ its tangent vector. We can take a side tangential vector $\sigma 
$ in $T_{x}M$ such that $\left\{ T,\sigma \right\} $ is a positive oriented
base. Therefore we have a decomposition:%
\begin{equation*}
r^{\prime \prime }=\frac{d^{2}r}{ds^{2}}=\kappa N=\kappa _{g}\sigma +\kappa
_{n}\xi ,
\end{equation*}%
where $N$, $\kappa _{g}$ and $\kappa _{n}$ are the normal vector, \textit{%
geodesic }and \textit{normal} \textit{curvatures} of $r$ on $M$,
respectively. Put $\sigma =a_{1}r_{u_{1}}+a_{2}r_{u_{2}}.$ Due to $%
T=r_{u_{1}}\frac{du_{1}}{ds}+u_{2}\frac{du_{2}}{ds}$ and $\left\langle
T,\sigma \right\rangle =0,$ we get%
\begin{equation*}
a_{1}=\theta \left( g_{12}\frac{du_{1}}{ds}+g_{22}\frac{du_{2}}{ds}\right) ,%
\text{ }a_{2}=-\theta \left( g_{11}\frac{du_{1}}{ds}+g_{12}\frac{du_{2}}{ds}%
\right) ,
\end{equation*}%
where $\theta =\theta \left( u_{1},u_{2}\right) $ is some nonzero smooth
function. Then we achieve%
\begin{equation*}
1=\det \left( \tilde{T},\tilde{\sigma}\right) =-\sqrt{\left\vert \det \left(
g_{ij}\right) \right\vert }\theta
\end{equation*}%
and hence%
\begin{equation*}
\sigma =-\frac{1}{\sqrt{\left\vert \det \left( g_{ij}\right) \right\vert }}%
\left[ \left( g_{12}\frac{du_{1}}{ds}+g_{22}\frac{du_{2}}{ds}\right)
r_{u_{1}}-\left( g_{11}\frac{du_{1}}{ds}+g_{12}\frac{du_{2}}{ds}\right)
r_{u_{2}}\right] .
\end{equation*}%
Accordingly, we compute that 
\begin{equation*}
\left. 
\begin{array}{l}
\kappa _{n}=\det \left( r^{\prime },\sigma ,r^{\prime \prime }\right) =%
\dfrac{1}{\sqrt{\left\vert \det \left( g_{ij}\right) \right\vert }}\det
\left( r_{u_{1}},r_{u_{2}},r^{\prime \prime }\right) \\ 
=\dfrac{1}{\sqrt{\left\vert \det \left( g_{ij}\right) \right\vert }}%
\sum\limits_{i,j=1}^{2}\det \left( r_{u_{1}},r_{u_{2}},r_{u_{i}u_{j}}\right)
\left( \dfrac{du_{i}}{ds}\right) \left( \dfrac{du_{j}}{ds}\right) ,%
\end{array}%
\right.
\end{equation*}%
which leads to the components of the second fundamental form given by 
\begin{equation}
h_{ij}=\frac{\det \left( r_{u_{1}},r_{u_{2}},r_{u_{i}u_{j}}\right) }{\sqrt{%
\left\vert \det \left( g_{ij}\right) \right\vert }},\text{ }r_{u_{i}u_{j}}=%
\frac{\partial ^{2}r}{\partial u_{i}\partial u_{j}},\text{ }i,j=1,2.  \notag
\end{equation}

Thus the \textit{Gaussian curvature} and the \textit{mean curvature} of $M$
are respectively defined by%
\begin{equation}
K=\frac{\det \left( h_{ij}\right) }{\det \left( g_{ij}\right) }  \tag{4.1}
\end{equation}%
and%
\begin{equation}
H= \frac{g_{11}h_{22}-2g_{12}h_{12}+g_{22}h_{11}}{2\det \left( g_{ij}\right) 
}.  \tag{4.2}
\end{equation}

By permutation of the coordinates, two different types of graph surfaces
appear up to the absolute of $\mathbb{I}_{p}^{3}.$ For a graph of the
function $u=u\left( x,y\right) ,$ the formulas (4.1) and (4.2) reduce to%
\begin{equation*}
K=-u_{xx}u_{yy}+\left( u_{xy}\right) ^{2},\text{ }H=\frac{1}{2}\left(
u_{xx}-u_{yy}\right).
\end{equation*}
Since the metric on the graph surface induced from $\mathbb{I}_{p}^{3}$ is $%
g=dx^{2}-dy^{2}$, it always becomes a flat surface. So, its Laplacian turns
to%
\begin{equation*}
\bigtriangleup =\frac{\partial ^{2}}{\partial u_{1}^{2}}-\frac{\partial ^{2}%
}{\partial u_{2}^{2}}.
\end{equation*}

On the other side, the Gaussian and mean curvatures of the graph of $%
u=u\left( y,z\right) $ are given by%
\begin{equation*}
K=-\frac{u_{yy}u_{zz}-\left( u_{yz}\right) ^{2}}{\left( u_{z}\right) ^{4}},%
\text{ }H=\frac{\left( u_{z}\right) ^{2}u_{yy}-2u_{y}u_{z}u_{yz}+\left(
\left( u_{y}\right) ^{2}-1\right) u_{zz}}{2\left( u_{z}\right) ^{3}}.
\end{equation*}

\section{Constant curvature surfaces of revolution in $\mathbb{I}_{p}^{3}$}

Da Silva \cite{8} provided via hyperbolic numbers that the pseudo isotropic
motion given by $\bar{x}=px,$ $\bar{y}=\frac{1}{p}y,$ $p\neq 0$ is
equivalent to the hyperpolic rotation (about $z-$axis) given by 
\begin{equation}
\bar{x}=x\cosh \theta +y\sinh \theta ,\text{ } \bar{y}=x\sinh \theta +y\cosh \theta , \tag{5.1}
\end{equation}
where $\theta \in \mathbb{R}.$ 

Let $u\longmapsto \left( u,0,f\left( u\right) \right) $ be a spacelike
admissible curve lying in the isotropic $xz-$plane of $\mathbb{I}_{p}^{3}$
for a smooth function $f$. Rotating it around $z-$axis via hyperbolic
rotations given by $\left( 5.1\right) $ we derive 
\begin{equation}
r\left( u,v\right) =\left( u\cosh v,u\sinh v,f\left( u\right) \right) . 
\tag{5.2}
\end{equation}%
We call the rotating curve \textit{profile curve. }If the profile curve is a
timelike curve $u\longmapsto \left( 0,u,f\left( u\right) \right) $ lying in
the isotropic $yz-$plane of $\mathbb{I}_{p}^{3}$, then rotating it around $%
\mathbf{z}-$axis yields 
\begin{equation}
r\left( u,v\right) =\left( u\sinh v,u\cosh v,f\left( u\right) \right) . 
\tag{5.3}
\end{equation}%
The surfaces given by $\left( 5.2\right) $ and $\left( 5.3\right) $ are
called \textit{surfaces of revolution} in $\mathbb{I}_{p}^{3}.$ The Gaussian
curvature of these surfaces in $\mathbb{I}_{p}^{3}$ is%
\begin{equation}
K=\frac{f^{\prime }f^{\prime \prime }}{u},  \tag{5.4}
\end{equation}%
where $f^{\prime }\left( u\right) =\dfrac{df}{du},$ etc.

Now we assume that it has nonzero constant Gaussian curvature\textit{\ }$%
K_{0}$ in $\mathbb{I}_{p}^{3}.$ Then $\left( 5.4\right) $ can be rewritten as%
\begin{equation}
f^{\prime }=\sqrt{c_{1}+K_{0}u^{2}},c_{1}\in \mathbb{R}.  \tag{5.5}
\end{equation}%
After integrating $\left( 5.5\right) ,$ we obtain%
\begin{equation*}
f\left( u\right) =\frac{u}{2}\sqrt{c_{1}+K_{0}u^{2}}+\frac{c_{1}}{2\sqrt{%
K_{0}}}\ln \left( 2K_{0}u+2\sqrt{K_{0}}\sqrt{c_{1}+K_{0}u^{2}}\right) +c_{2},%
\text{ }c_{2}\in \mathbb{R}
\end{equation*}%
which implies the following result.

\begin{theorem}
Let $M$ be a surface of revolution in $\mathbb{I}_{p}^{3}$ with nonzero
constant Gaussian curvature $K_{0}.$ Then its profile curve is of the form $%
\left( u,0,f\left( u\right) \right) $, where%
\begin{equation*}
f\left( u\right) =\frac{u}{2}\psi \left( u\right) +\frac{c_{1}}{2\sqrt{K_{0}}%
}\ln \left\vert 2\sqrt{K_{0}}\left( \sqrt{K_{0}}x+\psi \left( u\right)
\right) \right\vert
\end{equation*}%
for $\psi \left( u\right) =\sqrt{c_{1}+K_{0}u^{2}}$, $c_{1},c_{2}\in \mathbb{%
R}.$
\end{theorem}

We immediately have the following from $(5.4)$.

\begin{corollary}
A surface of revolution is flat in $\mathbb{I}_{p}^{3}$ if and only if its
profile curve is a non-isotropic line given by $\left(
u,0,c_{1}u+c_{2}\right) ,$ $c_{1},c_{2}\in \mathbb{R}.$
\end{corollary}

The mean curvature $H$ of a surface of revolution $M$ in $\mathbb{I}_{p}^{3} 
$ is%
\begin{equation}
H=\frac{1}{2}\left( \frac{f^{\prime }}{u}+f^{\prime \prime }\right) . 
\tag{5.6}
\end{equation}%
Assume that $M$ has constant mean curvature $H_{0}.$ After solving $\left(
5.6\right) $ we deduce%
\begin{equation*}
f\left( u\right) =\frac{H_{0}}{2}u^{2}+c_{1}\ln u+c_{2},\text{ }%
c_{1},c_{2}\in \mathbb{R}.
\end{equation*}

Therefore we have proved the following results.

\begin{theorem}
Let $M$ be a surface of revolution in $\mathbb{I}_{p}^{3}$ with constant
mean curvature $H_{0}.$ Then its profile curve is of the form $\left(
u,0,f\left( u\right) \right) $, where%
\begin{equation*}
f\left( u\right) =\frac{H_{0}}{2}u^{2}+c_{1}\ln u+c_{2},\text{ }%
c_{1},c_{2}\in \mathbb{R}.
\end{equation*}
\end{theorem}

\begin{corollary}
A surface of revolution is minimal in $\mathbb{I}_{p}^{3}$ if and only if
its profile curve is a non-isotropic curve given by $\left( u,0,c_{1}\ln
u+c_{2}\right) ,$ $c_{1},c_{2}\in \mathbb{R}.$
\end{corollary}

\begin{example}
Take the surfaces of revolution in $\mathbb{I}_{p}^{3}$ parameterized%
\begin{equation*}
r\left( u,v\right) =\left( u\cosh v,u\sinh v,u\right) ,\text{ }\left(
u,v\right) \in \left[ 1,2\right] \times \left[ 0,1\right]
\end{equation*}%
and 
\begin{equation*}
\left( u,v\right) =\left( u\cosh v,u\sinh v,\ln u+u^{2}\right) ,\text{ }%
\left( u,v\right) \in \left[ 1,2\right] \times \left[ -1,1\right] .
\end{equation*}%
The above first surface is flat and the second is a constant mean curvature
surface of revolution, $H=2.$ We plot these as in Figure 2 and Figure 3,
respectively.
\end{example}

\begin{figure}[ht]
\begin{center}
\includegraphics[scale=0.5]{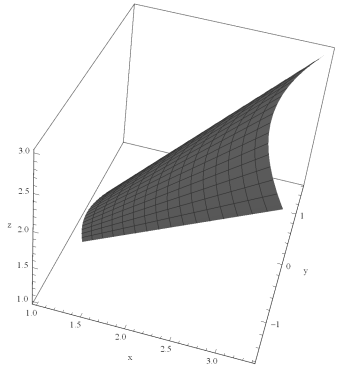}
\end{center}
\caption{A flat surface of revolution, $K=0$.}
\end{figure}

\begin{figure}[ht]
\begin{center}
\includegraphics[scale=0.5]{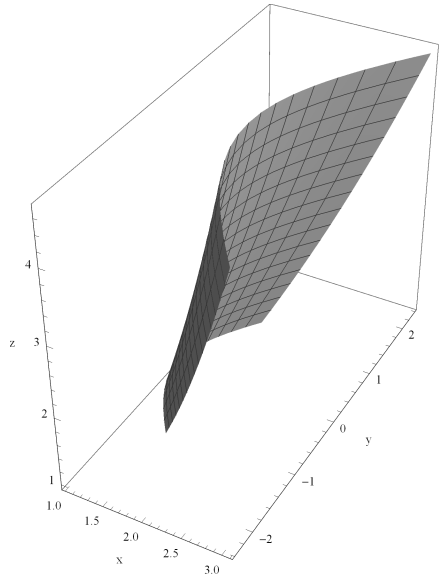}
\end{center}
\caption{A constant curvature surface of revolution, $H=2$.}
\end{figure}

\section{Surfaces of revolution with $H^{2}=K$ in $\mathbb{I}_{p}^{3}$}

Next we aim to classify the surfaces of revolution given by (5.2) in $%
\mathbb{I}_{p}^{3}$ that satisfy $H^{2}=K$ which is the equality sign of the Euler inequality. For more generalizations of the famous inequality, see \cite{4,19,20}.

By considering the equalities $\left( 5.4\right) $ and $\left( 5.6\right) ,$ we have%
\begin{equation}
\frac{1}{4}\left( \left( \frac{f^{\prime }}{u}\right) ^{2}+2\frac{f^{\prime
}f^{\prime \prime }}{u}+\left( f^{\prime \prime }\right) ^{2}\right) =\frac{%
f^{\prime }f^{\prime \prime }}{u}.  \tag{6.1}
\end{equation}%
We can rewirte $\left( 6.1\right) $ as%
\begin{equation*}
\left( \frac{f^{\prime }}{u}-f^{\prime \prime }\right) ^{2}=0,
\end{equation*}%
which implies%
\begin{equation*}
\frac{f^{\prime }}{u}-f^{\prime \prime }=0.
\end{equation*}%
After solving this, we obtain%
\begin{equation*}
f\left( u\right) =c_{1}\frac{u^{2}}{2}+c_{2}
\end{equation*}%
for $c_{1},c_{2}\in 
\mathbb{R}
.$ By comparing (5.2) with (6.2) we see that the surface of revolution can be given in explicit form

\begin{equation}
z=\frac{c_{1}}{2}\left( x^{2}-y^{2}\right) +c_{2},  \tag{6.3}
\end{equation}%
which implies the following result.

\begin{theorem}
The surfaces of revolution given by (5.2) in $\mathbb{I}_{p}^{3}$ with $%
H^{2}=K$ are only the spheres of parabolic type.
\end{theorem}

\begin{example}
Consider the sphere of parabolic type in $\mathbb{I}_{p}^{3}$ given via $%
\left( 6.3\right) $ such that $c_{1}=2$ and $c_{2}=0.$ Then its curvatures
become $H=2$ and $K=4$. We plot it as in Figure 4.
\end{example}

\begin{figure}[th]
\begin{center}
\includegraphics[scale=0.5]{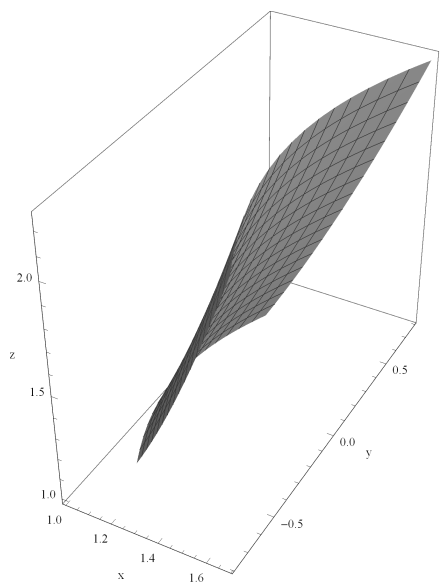}
\end{center}
\caption{A surface of revolution with $H^{2}=K$.}
\end{figure}

\section{Acknowlodgements}

The figures in the present study were made by Wolfram Mathematica 11.0.

\end{document}